\documentclass[11pt]{article}
\usepackage{amsfonts}
\usepackage{graphics}
\usepackage{indentfirst}
\usepackage{color}
\usepackage{cite}
\usepackage{latexsym}
\usepackage[paper=a4paper, left=2.1cm, right=2.1cm, top=2.2cm, bottom=1.5cm, headheight=5.5pt, footskip=0.8cm, footnotesep=0.8cm, centering, includefoot]{geometry}
\usepackage{amsmath}
\allowdisplaybreaks
\usepackage{amssymb}
\usepackage[dvips]{epsfig}
\usepackage{amscd}

\newtheorem{theorem}{Theorem}[section]
\newtheorem{remark}{Remark}[section]
\newtheorem{definition}{Definition}[section]
\newtheorem{lemma}{Lemma}[section]

\DeclareMathOperator{\curl}{curl}
\DeclareMathOperator{\loc}{loc}

\DeclareMathOperator{\divv}{div}

\makeatletter
\@addtoreset{equation}{section}
\makeatother
\makeatletter
\@addtoreset{equation}{section}
\makeatother

\title{{Singularity formation to the 2D Cauchy problem of the full compressible Navier-Stokes equations with zero heat conduction}
\thanks{Supported by China Postdoctoral Science Foundation (No. 2017M610579), Fundamental Research Funds for the Central Universities (No. XDJK2017C050) and the Doctoral Fund of Southwest University (No. SWU116033).}
}

\author{Xin Zhong\thanks{School of Mathematics and Statistics, Southwest University, Chongqing 400715,
People's Republic of China ({\tt xzhong1014@amss.ac.cn}).
}
}
\date{ }

\begin{document}
\maketitle

\begin{abstract}
The formation of singularity and breakdown of strong solutions to the two-dimensional (2D) Cauchy problem of the full compressible Navier-Stokes equations with zero heat conduction are considered. It is shown that for the initial density allowing vacuum, the strong solution exists globally if the density $\rho$ and the pressure $P$ satisfy $\|\rho\|_{L^{\infty}(0,T;L^\infty)}+\|P\|_{L^{\infty}(0,T;L^\infty)}<\infty$.
In addition, the initial density can even have compact support. The logarithm-type estimate for the Lam{\'e} system and some weighted estimates play a crucial role in the proof.
\end{abstract}


Math Subject Classification: 35Q30; 35B65

\section{Introduction}
Let $\Omega\subset\mathbb{R}^2$ be a domain, the motion of a viscous, compressible, and heat conducting Navier-Stokes flow in $\Omega$ can be described by the full compressible Navier-Stokes equations
\begin{align}\label{1.1}
\begin{cases}
\rho_{t}+\divv(\rho\mathbf{u})=0,\\
(\rho\mathbf{u})_{t}+\divv(\rho\mathbf{u}\otimes\mathbf{u})
-\mu\Delta\mathbf{u}
-(\lambda+\mu)\nabla\divv\mathbf{u}+\nabla P=
\mathbf{0},\\
c_{\nu}[(\rho\theta)_{t}+\divv(\rho\mathbf{u}\theta)]
+P\divv\mathbf{u}-\kappa\Delta\theta
=2\mu|\mathfrak{D}(\mathbf{u})|^2+\lambda(\divv\mathbf{u})^2.
\end{cases}
\end{align}
Here, $t\geq0$ is the time, $x\in\Omega$ is the spatial coordinate, and the unknown functions $\rho, \mathbf{u}, P=R\rho\theta\ (R>0), \theta$ are the fluid density, velocity, pressure, and the absolute temperature respectively; $\mathfrak{D}(\mathbf{u})$ denotes the deformation tensor given by
\begin{equation*}
\mathfrak{D}(\mathbf{u})=\frac{1}{2}(\nabla\mathbf{u}+(\nabla\mathbf{u})^{tr}).
\end{equation*}
The constant viscosity coefficients $\mu$ and $\lambda$ satisfy the physical restrictions
\begin{equation}\label{1.2}
\mu>0,\ \mu+\lambda\geq0.
\end{equation}
Positive constants $c_\nu$ and $\kappa$ are respectively the heat capacity, the ratio of the heat conductivity coefficient over the heat capacity.

There is huge literature on the studies about the theory of
well-posedness of solutions to the Cauchy problem and the initial boundary value problem (IBVP) for the compressible Navier-Stokes system due to the physical importance, complexity, rich phenomena and mathematical challenges, refer to \cite{CK2006,CK20062,H1995,H1997,HLX2012,HL2013,L1998,F2004,NS2004,
LX2013,LL2014,YZ2017,WZ2017,W2017,MN1980,MN1983} and references therein.
In particular, non-vacuum small perturbations of a uniform non-vacuum
constant state have been shown existing globally in time and remain smooth in any space dimensions \cite{MN1980,MN1983}, while for general data which may contain vacuum states, only weak solutions are shown to exist for the compressible Navier-Stokes system in multi-dimension with special equation of state as in \cite{F2004,FNP2001,L1998}, yet the uniqueness
and regularity of these weak solutions remain unknown.
Despite the surprising results on global well-posedness of the strong (or classical) solution to the multi-dimensional compressible Navier-Stokes system for initial data with small total energy but possible large oscillations and containing vacuum states \cite{HLX2012,HL2013,LX2013,YZ2017,WZ2017},
it is an outstanding challenging open problem to investigate the global well-posedness for general large strong solutions with vacuum.

%However, many physical important and mathematical fundamental problems are still open due to the lack of smoothing mechanism and the strong nonlinearity. When the initial density allows vacuum, the local large strong solutions to the Cauchy problem of 3D full Navier-Stokes equations have been obtained by Cho-Kim \cite{CK2006}.
%For the global well-posedness of strong (or classical) solutions, Huang-Li-Xin \cite{HLX2012} and Li-Xin \cite{LX2013} established the global existence and uniqueness of strong (or classical) solutions to the 3D and 2D isentropic Navier-Stokes equations, respectively, provided the smooth initial data are of small total energy. Furthermore, the authors \cite{L1998,FNP2001} showed the global existence of renormalized solutions to the compressible Navier-Stokes equations for general large initial data.

Therefore, it is important to study the mechanism of blow-up and structure of possible singularities of strong (or classical) solutions to the compressible Navier-Stokes equations. The pioneering work can be traced
to Serrin's criterion \cite{S1962} on the Leray-Hopf weak solutions to the three-dimensional incompressible Navier-Stokes equations, which can be stated
that if a weak solution $\mathbf{u}$ satisfies
\begin{equation}\label{1.3}
\mathbf{u}\in L^s(0,T;L^r),\ \text{for}\ \frac2s+\frac3r=1,\ 3<r\leq\infty,
\end{equation}
then it is regular.

Recently, there are several results on the blow-up criteria of strong (or classical) solutions to the compressible Navier-Stokes equations. Precisely, let $0<T^*<+\infty$ be the maximum time of existence of strong solutions.
For the 3D isentropic flows, Huang--Li-Xin \cite{HLX20112} obtained the following Serrin type criterion
\begin{equation}\label{1.4}
\lim_{T\rightarrow T^{*}}\left(\|\divv\mathbf{u}\|_{L^1(0,T;L^\infty)}+\|\mathbf{u}\|_{L^s(0,T;L^r)}\right)=\infty,
\end{equation}
where $r$  and $s$ as in \eqref{1.3}. In \cite{HLX20111}, they also proved a Beale-Kato-Majda type criterion as follows
\begin{equation}\label{1.04}
\lim_{T\rightarrow T^{*}}\|\mathfrak{D}(\mathbf{u})\|_{L^1(0,T;L^\infty)}=\infty.
\end{equation}
For the IBVP of 3D isentropic flows, Sun-Wang-Zhang \cite{SWZ20111} established
\begin{equation}
\lim_{T\rightarrow T^{*}}\|\rho\|_{L^\infty(0,T;L^\infty)}=\infty
\end{equation}
provided that
\begin{equation}\label{1.6}
7\mu>\lambda.
\end{equation}
For the 3D full compressible Navier-Stokes equations, under the condition \eqref{1.6}, Fan-Jiang-Ou \cite{FJO2010} showed that
\begin{equation}\label{1.5}
\lim_{T\rightarrow T^{*}}\left(\|\nabla\mathbf{u}\|_{L^1(0,T;L^\infty)}+\|\theta\|_{L^\infty(0,T;L^\infty)}\right)
=\infty.
\end{equation}
Under just the physical condition
\begin{equation}\label{1.06}
\mu>0,\ 2\mu+3\lambda\geq0,
\end{equation}
Huang-Li-Wang established the criterion \eqref{1.4} for the 3D barotropic case \cite{HLX20112} still holds for the full Navier-Stokes system.
For the Cauchy problem and the IBVP of 3D full compressible Navier-Stokes system, Huang-Li \cite{HL2013} proved that
\begin{equation}\label{1.7}
\lim_{T\rightarrow T^{*}}\left(\|\rho\|_{L^\infty(0,T;L^\infty)}+\|\mathbf{u}\|_{L^s(0,T;L^r)}\right)
=\infty,\ \text{for}\ \frac2s+\frac3r\leq1,\ 3<r\leq\infty.
\end{equation}
However, for the IBVP of 2D full Navier-Stokes equations, Wang \cite{W2014} showed the formation of singularity must be caused by losing the bound of $\divv\mathbf{u}$. More precisely, she obtained
\begin{equation}\label{1.08}
\lim_{T\rightarrow T^{*}}\|\divv\mathbf{u}\|_{L^1(0,T;L^\infty)}=\infty.
\end{equation}
For more information on the blow-up criteria of compressible flows, we refer to \cite{CJ2006,R2008,SWZ20112,WZ2013,X1998,XY2013,HX2016} and the references therein.

It is worth noting that one would not expect better regularities of the solutions of
\eqref{1.1} in general because of Xin's result \cite{X1998}, where the author proved that there is no global smooth solution to the Cauchy problem of \eqref{1.1} if the initial density is nontrivial compactly supported. Very recently, Liang-Shi \cite{LS2015} obtained the local existence of strong (or classical) solutions for the non-isentropic compressible Navier-Stokes equations without heat-conductivity.
These motivate us to find some possible blow-up criterion of regular solutions to the system \eqref{1.1} with zero heat conduction, especially of strong solutions.
%a natural question arises: how we can describe the mechanism of blow up to ensure the global existence of strong (or classical) solutions to the system \eqref{1.1} with zero heat conduction?
%However, the mechanism of blow-up of strong solutions to the problem \eqref{1.1} with $\kappa=0$ is still unknown.
In fact, this is the main aim of this paper.

When $\kappa=0$, and without loss of generality, take $c_\nu=R=1$, the system \eqref{1.1} can be written as
\begin{align}\label{1.10}
\begin{cases}
\rho_{t}+\divv(\rho\mathbf{u})=0,\\
(\rho\mathbf{u})_{t}+\divv(\rho\mathbf{u}\otimes\mathbf{u})
-\mu\Delta\mathbf{u}
-(\lambda+\mu)\nabla\divv\mathbf{u}+\nabla P=
\mathbf{0},\\
P_{t}+\divv(P\mathbf{u})
+P\divv\mathbf{u}
=2\mu|\mathfrak{D}(\mathbf{u})|^2+\lambda(\divv\mathbf{u})^2.
\end{cases}
\end{align}
The present paper is aimed at giving a blow-up criterion of strong solutions to the Cauchy problem  of the system \eqref{1.10} with the initial condition
\begin{equation}\label{1.11}
(\rho,\mathbf{u},P)(x,0)=(\rho_0,\mathbf{u}_0,P_0)(x),\ \ x\in\mathbb{R}^2,
\end{equation}
and the far field behavior
\begin{equation}\label{1.12}
(\rho,\mathbf{u},P)(x,t)\rightarrow(0, \mathbf{0},0),\ \text{as}\ |x|\rightarrow+\infty,\ t>0.
\end{equation}

Before stating our main result, we first explain the notations and conventions used throughout this paper. For $r>0$, set
\begin{equation*}
B_r \triangleq\left.\left\{x\in\mathbb{R}^2\right|\,|x|<r \right\},
\quad \int \cdot dx\triangleq\int_{\mathbb{R}^2}\cdot dx.
\end{equation*}
For $1\leq p\leq\infty$ and integer $k\geq0$, the standard Sobolev spaces are denoted by:
\begin{equation*}
L^p=L^p(\mathbb{R}^2),\ W^{k,p}=W^{k,p}(\mathbb{R}^2), \ H^{k}=H^{k,2}(\mathbb{R}^2),
\ D^{k,p}=\{u\in L_{\loc}^1|\nabla^k u\in L^p\}.
\end{equation*}
Now we define precisely what we mean by strong solutions to the problem \eqref{1.10}--\eqref{1.12}.
\begin{definition}[Strong solutions]\label{def1}
$(\rho,\mathbf{u},P)$ is called a strong solution to \eqref{1.10}--\eqref{1.12} in $\mathbb{R}^2\times(0,T)$, if for some $q_0>2$,
\begin{equation*}
\begin{split}
\begin{cases}
\rho\geq0,\ \rho\in C([0,T];L^1\cap H^1\cap W^{1,q_0}),\ \rho_t\in C([0,T];L^{q_0}),\\
\mathbf{u}\in C([0,T];D^{1,2}\cap D^{2,2})\cap L^{2}(0,T;D^{2,q_0}), \\
\mathbf{u}_t\in L^{2}(0,T;D^{1,2}),\
\sqrt{\rho}\mathbf{u}_{t}\in L^{\infty}(0,T;L^{2}), \\
P\geq0,\ P\in C([0,T];L^1\cap H^1\cap W^{1,q_0}),\ P_t\in C([0,T];L^{q_0}), \\
\end{cases}
\end{split}
\end{equation*}
and $(\rho,\mathbf{u},P)$ satisfies both \eqref{1.10} almost everywhere in $\mathbb{R}^2\times(0,T)$ and \eqref{1.11} almost everywhere in $\mathbb{R}^2$.
\end{definition}

Without loss of generality, we assume that the initial density $\rho_0$ satisfies
\begin{equation}\label{1.8}
\int_{\mathbb{R}^2} \rho_0dx=1,
\end{equation}
which implies that there exists a positive constant $N_0$ such that
\begin{equation}\label{1.9}
\int_{B_{N_0}}  \rho_0  dx\ge \frac12\int\rho_0dx=\frac12.
\end{equation}

Our main result reads as follows:
\begin{theorem}\label{thm1.1}
Let $\eta_0$ be a positive number and
\begin{equation}\label{2.01}
\bar{x}\triangleq(e+|x|^2)^{1/2}\log^{1+\eta_0}(e+|x|^2).
\end{equation}
In addition to \eqref{1.8} and \eqref{1.9}, assume that the initial data $(\rho_0\geq0, \mathbf{u}_0,P_0\geq0)$ satisfies for any given numbers $a>1$ and $q>2$,
\begin{equation}\label{2.02}
\rho_{0}\bar{x}^{a}\in L^{1}\cap H^{1}\cap W^{1,q},\ \nabla\mathbf{u}_{0}\in H^1,\  \sqrt{\rho_0}\mathbf{u}_0\in L^2, \ P_0\in L^{1}\cap H^{1}\cap W^{1,q},
\end{equation}
and the compatibility conditions
\begin{equation}\label{A2}
-\mu\Delta\mathbf{u}_0-(\lambda+\mu)\nabla\divv\mathbf{u}_0+\nabla P_0=\sqrt{\rho_0}\mathbf{g}
\end{equation}
for some $\mathbf{g}\in L^2(\Omega)$.
Let $(\rho,\mathbf{u},P)$ be a strong solution to the problem \eqref{1.10}--\eqref{1.12}. If $T^{*}<\infty$ is the maximal time of existence for that solution,
then we have
\begin{align}\label{B}
\lim_{T\rightarrow T^{*}}\left(\|\rho\|_{L^{\infty}(0,T;L^\infty)}
+\|P\|_{L^{\infty}(0,T;L^\infty)}\right)=\infty.
\end{align}
\end{theorem}

Several remarks are in order.
\begin{remark}\label{re1.1}
The local existence of a strong solution with initial data as in Theorem \ref{thm1.1} was established in \cite{LS2015,LH2016}. Hence, the maximal time $T^{*}$ is well-defined.
\end{remark}
\begin{remark}\label{re1.2}
According to \eqref{B}, the upper bound of the temperature $\theta$ is not the key point to make sure that the solution $(\rho,\mathbf{u},P)$ is a global one, and it may go to infinity in the vacuum region within the life span of our strong solution.
\end{remark}
\begin{remark}\label{re1.3}
Compared with \cite{HX2016}, where the authors investigated blow-up criteria for the 3D Cauchy problem and the IBVP of non-isentropic Navier-Stokes equations with zero heat conduction, there is no need to impose additional restrictions on the viscosity coefficients $\mu$ and $\lambda$ except the physical restrictions \eqref{1.2}.
\end{remark}

We now make some comments on the analysis of this paper. We mainly make use of continuation argument to prove Theorem \ref{thm1.1}. That is, suppose that \eqref{B} were false, i.e.,
\begin{equation*}
\lim_{T\rightarrow T^*}\left(\|\rho\|_{L^{\infty}(0,T;L^\infty)}
+\|P\|_{L^{\infty}(0,T;L^\infty)}\right)\leq M_0<\infty.
\end{equation*}
We want to show that
\begin{equation*}
\sup_{0\leq t\leq T^*}\left(\|(\rho,P)\|_{H^1\cap W^{1,q}}
+\|\rho\bar{x}^{a}\|_{L^{1}\cap H^{1}\cap W^{1,q}}
+\|\nabla\mathbf{u}\|_{H^1}\right) \leq C<+\infty.
\end{equation*}

It should be pointed out that the crucial techniques of proofs in \cite{W2014} cannot be adapted directly to the situation treated here, since their arguments depend crucially on the boundedness of the domains and $\kappa>0$. Moreover, technically, it is hard to modify the three-dimensional analysis of \cite{HX2016} to the two-dimensional case
with initial density containing vacuum since the
analysis of \cite{HX2016} depends crucially on the a priori $L^6$-bound on the velocity, while in two dimensions it seems difficult to bound the $L^p(\mathbb{R}^2)$-norm of $\mathbf{u}$ just in terms of $\|\sqrt{\rho}\mathbf{u}\|_{L^{2}(\mathbb{R}^2)}$ and $\|\nabla\mathbf{u}\|_{L^{2}(\mathbb{R}^2)}$ for any $p\geq1$.
To overcome these difficulties mentioned above, some new ideas are needed.
Inspired by \cite{LX2013,LSZ2015}, we first observe that if the initial density decays not too slow at infinity, i.e., $\rho_0\bar{x}^a\in L^1(\mathbb{R}^2)$ for some positive constant $a>1$ (see \eqref{2.02}), then for any $\eta\in(0,1]$, we can show that (see \eqref{06.2})
\begin{equation}\label{1.20}
\mathbf{u}\bar{x}^{-\eta}\in L^{p_0}(\mathbb{R}^2),\  \text{for some}\  p_0>1.
\end{equation}
Then, motivated by \cite{H1995,HX2016}, in order to get the $L_t^\infty L_x^{2}$-norm of $\sqrt{\rho}\dot{\mathbf{u}}$, we first show the desired a priori estimates of the $L_t^\infty L_x^2$-norm of $\nabla\mathbf{u}$, which is the second key observation in this paper (see Lemma \ref{lem35}). Next, the a priori estimates on the $L_t^\infty L_x^{q}$-norm of $(\nabla\rho,\nabla P)$ and the $L_t^1 L_x^{\infty}$-norm of the velocity gradient can be obtained (see Lemma \ref{lem37}) simultaneously by solving a logarithm Gronwall inequality based on a logarithm estimate for the Lam{\'e} system (see Lemma \ref{lem25}). Finally, with the help of \eqref{1.20}, we can get the spatial weighted estimate of the density (see Lemma \ref{lem38}).

The rest of this paper is organized as follows. In Section \ref{sec2}, we collect some elementary facts and inequalities that will be used later. Section \ref{sec3} is devoted to the proof of Theorem \ref{thm1.1}.

\section{Preliminaries}\label{sec2}

In this section, we will recall some known facts and elementary inequalities that will be used frequently later.

We begin with the following Gronwall's inequality (see \cite[pp. 12--13]{T2006}), which plays a central role in proving a priori estimates on strong solutions $(\rho,\mathbf{u},P)$.
\begin{lemma}\label{lem21}
Suppose that $h$ and $r$ are integrable on $(a, b)$ and nonnegative a.e. in $(a, b)$. Further assume that $y\in C[a, b], y'\in L^1(a, b)$, and
\begin{equation*}
y'(t)\leq h(t)+r(t)y(t)\ \ \text{for}\ a.e\ t\in(a,b).
\end{equation*}
Then
\begin{equation*}
y(t)\leq \left[y(a)+\int_{a}^{t}h(s)\exp\left(-\int_{a}^{s}r(\tau)d\tau\right)ds\right]
\exp\left(\int_{a}^{t}r(s)ds\right),\ \ t\in[a,b].
\end{equation*}
\end{lemma}

Next, the following Gagliardo-Nirenberg inequality (see \cite{N1959}) will be used later.
\begin{lemma}[Gagliardo-Nirenberg]\label{lem22}
For $p\in[2,\infty), r\in(2,\infty)$, and $s\in(1,\infty)$, there exists some generic constant $C>0$ which may depend on $p,$ $r$, and $s$ such that for $f\in H^{1}(\mathbb{R}^2)$ and $g\in L^{s}(\mathbb{R}^2)\cap D^{1,r}(\mathbb{R}^2)$, we have
\begin{eqnarray*}
& & \|f\|_{L^p(\mathbb{R}^2)}^{p}\leq C\|f\|_{L^2(\mathbb{R}^2)}^{2}\|\nabla f\|_{L^2(\mathbb{R}^2)}^{p-2}, \\
& & \|g\|_{C(\overline{\mathbb{R}^2})}\leq C\|g\|_{L^s(\mathbb{R}^2)}^{s(r-2)/(2r+s(r-2))}\|\nabla g\|_{L^r(\mathbb{R}^2)}^{2r/(2r+s(r-2))}.
\end{eqnarray*}
\end{lemma}

The following weighted $L^m$ bounds for elements of the Hilbert space $\tilde{D}^{1,2}(\mathbb{R}^2)\triangleq\{ v \in H_{\loc}^{1}(\mathbb{R}^2)|\nabla v \in L^{2}(\mathbb{R}^2)\}$ can be found in \cite[Theorem B.1]{L1996}.
\begin{lemma}\label{1leo}
For $m\in [2,\infty)$ and $\theta\in (1+m/2,\infty),$ there exists a positive constant $C$ such that for all $v\in  \tilde{D}^{1,2}(\mathbb{R}^2),$
\begin{equation}\label{3h}
\left(\int_{\mathbb{R}^2} \frac{|v|^m}{e+|x|^2}\left(\log \left(e+|x|^2\right)\right)^{-\theta}dx  \right)^{1/m}
\leq C\|v\|_{L^2(B_1)}+C\|\nabla v\|_{L^2(\mathbb{R}^2) }.
\end{equation}
\end{lemma}

The combination of Lemma \ref{1leo} and the Poincar\'e inequality yields
the following useful results on weighted bounds, whose proof can be found in \cite[Lemma 2.4]{LX2013}.

\begin{lemma}\label{lem26}
Let $\bar x$ be as in \eqref{2.01}. Assume that $\rho \in L^1(\mathbb{R}^2)\cap L^\infty(\mathbb{R}^2)$ is a non-negative function such that
\begin{equation*} \label{2.i2}
\|\rho\|_{L^1(B_{N_1})} \ge M_1, \quad \|\rho\|_{L^1(\mathbb{R}^2)\cap L^\infty(\mathbb{R}^2)}\le M_2,
\end{equation*}
for positive constants $M_1, M_2$, and $ N_1\ge 1$. Then for $\varepsilon> 0$ and $\eta>0,$ there is a positive constant $C$ depending only on $\varepsilon,\eta, M_1,M_2$, and $ N_1$, such that every $v\in \tilde{D}^{1,2}(\mathbb{R}^2)$ satisfies
\begin{equation}\label{22}
\|v\bar x^{-\eta}\|_{L^{(2+\varepsilon)/\tilde{\eta}}(\mathbb{R}^2)}
\le C \|{\sqrt\rho}v\|_{L^2(\mathbb{R}^2)}+C \|\nabla v\|_{L^2(\mathbb{R}^2)},
\end{equation}
with $\tilde{\eta}=\min\{1,\eta\}$.
\end{lemma}

Next, the following Beale-Kato-Majda type inequality (see \cite[Lemma 2.3]{HLX20112}) will be used to estimate $\|\nabla\mathbf{u}\|_{L^\infty}$.
\begin{lemma}\label{lem25}
For $q\in(2,\infty)$, there is a constant $C(q)>0$ such that for all $\nabla\mathbf{v}\in L^2\cap D^{1,q}$, it holds that
\begin{equation}\label{2.2}
\|\nabla\mathbf{v}\|_{L^\infty}\leq C\left(\|\divv\mathbf{v}\|_{L^\infty}
+\|\curl\mathbf{v}\|_{L^\infty}\right)\log(e+\|\nabla^2\mathbf{v}\|_{L^q})
+C\|\nabla\mathbf{v}\|_{L^2}+C.
\end{equation}
\end{lemma}

Finally, for $\nabla^{\bot}\triangleq(-\partial_2,\partial_1)$, denoting the material derivative of $f$ by $\dot{f}\triangleq f_t+\mathbf{u}\cdot\nabla f$,
then we have the following $L^p$-estimate (see \cite[Lemma 2.5]{LX2013}) for the elliptic system derived from the momentum equations \eqref{1.10}$_2$:
\begin{equation}\label{3.014}
\Delta F=\divv(\rho\dot{\mathbf{u}}),\
\mu\Delta\omega=\nabla^{\bot}\cdot(\rho\dot{\mathbf{u}}),
\end{equation}
where $F$ is the effective viscous flux, $\omega$ is vorticity given by
\begin{equation}\label{3.14}
F=(\lambda+2\mu)\divv\mathbf{u}-P,\
\omega=\partial_1u_2-\partial_2u_1.
\end{equation}
\begin{lemma}\label{lem23}
Let $(\rho, \mathbf{u}, P)$ be a smooth solution of \eqref{1.10}. Then for $p\geq2$ there exists a positive constant $C$ depending only on $p, \mu$ and $\lambda$ such that
\begin{align}
&\|\nabla F\|_{L^p}+\|\nabla\omega\|_{L^p}
\leq C\|\rho\dot{\mathbf{u}}\|_{L^p}, \label{2.3} \\
&\|F\|_{L^p}+\|\omega\|_{L^p}
\leq C\|\rho\dot{\mathbf{u}}\|_{L^2}^{1-\frac2p}
\left(\|\nabla\mathbf{u}\|_{L^2}+\|P\|_{L^2}\right)^{\frac2p}, \label{2.4} \\
&\|\nabla\mathbf{u}\|_{L^p}
\leq C\|\rho\dot{\mathbf{u}}\|_{L^2}^{1-\frac2p}
\left(\|\nabla\mathbf{u}\|_{L^2}+\|P\|_{L^2}\right)^{\frac2p}
+C\|P\|_{L^p}. \label{2.5}
\end{align}
\end{lemma}

\section{Proof of Theorem \ref{thm1.1}}\label{sec3}

Let $(\rho,\mathbf{u},P)$ be a strong solution described in Theorem \ref{thm1.1}. Suppose that \eqref{B} were false, that is, there exists a constant $M_0>0$ such that
\begin{equation}\label{3.1}
\lim_{T\rightarrow T^*}\left(\|\rho\|_{L^{\infty}(0,T;L^\infty)}
+\|P\|_{L^{\infty}(0,T;L^\infty)}\right)\leq M_0<\infty.
\end{equation}

First, the estimate on the $L^\infty(0,T;L^p)$-norm of the density could be deduced directly from \eqref{1.10}$_1$ and \eqref{3.1}.
\begin{lemma}\label{lem31}
Under the condition \eqref{3.1}, it holds that for any $T\in[0,T^*)$,
\begin{equation}\label{3.2}
\sup_{0\leq t\leq T}\|\rho\|_{L^1\cap L^\infty}\leq C,
\end{equation}
where and in what follows, $C,C_1,C_2$ stand for generic positive constants depending only on $M_0,\lambda,\mu,T^{*}$, and the initial data.
\end{lemma}

Next, we have the following estimate which is similar to the energy estimate.
\begin{lemma}\label{lem32}
Under the condition \eqref{3.1}, it holds that for any $T\in[0,T^*)$,
\begin{equation}\label{3.3}
\sup_{0\leq t\leq T}\left(\|\sqrt{\rho}\mathbf{u}\|_{L^2}^2+\|P\|_{L^1\cap L^\infty}
\right)+\int_{0}^T\|\nabla\mathbf{u}\|_{L^2}^2 dt \leq C.
\end{equation}
\end{lemma}
{\it Proof.}
It follows from \eqref{1.10}$_3$ that
\begin{equation}\label{7.01}
P_t+\mathbf{u}\cdot\nabla P+2P\divv\mathbf{u}
=F\triangleq2\mu|\mathfrak{D}(\mathbf{u})|^2
+\lambda(\divv\mathbf{u})^2\geq0.
\end{equation}
Define particle path before blowup time
\begin{align*}
\begin{cases}
\frac{d}{dt}\mathbf{X}(x,t)
=\mathbf{u}(\mathbf{X}(x,t),t),\\
\mathbf{X}(x,0)=x.
\end{cases}
\end{align*}
Thus, along particle path, we obtain from \eqref{7.01} that
\begin{align*}
\frac{d}{dt}P(\mathbf{X}(x,t),t)
=-2P\divv\mathbf{u}+F,
\end{align*}
which implies
\begin{equation}\label{3.03}
P(\mathbf{X}(x,t),t)=\exp\left(-2\int_{0}^t\divv\mathbf{u}ds\right)
\left[P_0+\int_{0}^t\exp\left(2\int_{0}^s\divv\mathbf{u}d\tau\right)Fds\right]\geq0.
\end{equation}

Next, multiplying \eqref{1.10}$_2$ by $\mathbf{u}$ and integrating over $\mathbb{R}^2$, we obtain after integrating by parts that
\begin{align}\label{3.4}
\frac{1}{2}\frac{d}{dt}\int\rho|\mathbf{u}|^2dx
+\int\left[\mu|\nabla\mathbf{u}|^2+(\lambda+\mu)(\divv\mathbf{u})^2
\right]dx=\int P\divv\mathbf{u}dx.
\end{align}
Integrating \eqref{1.10}$_3$ with respect to $x$ and then adding the resulting equality to \eqref{3.4} give rise to
\begin{align}\label{3.004}
\frac{d}{dt}\int\left(\frac12\rho|\mathbf{u}|^2+P\right)dx=0,
\end{align}
which combined with \eqref{3.03}, \eqref{2.02}, and \eqref{3.1} leads to
\begin{align}\label{3.04}
\sup_{0\leq t\leq T}\left(\|\sqrt{\rho}\mathbf{u}\|_{L^2}^2+\|P\|_{L^1\cap L^\infty}
\right)\leq C.
\end{align}
This together with \eqref{3.4} and Cauchy-Schwarz inequality yields
\begin{align}\label{3.003}
\frac{d}{dt}\|\sqrt{\rho}\mathbf{u}\|_{L^2}^2
+\mu\|\nabla\mathbf{u}\|_{L^2}^2
\leq C.
\end{align}
So the desired \eqref{3.3} follows from \eqref{3.04} and \eqref{3.003} integrated with respect to $t$. This completes the proof of Lemma \ref{lem32}.
\hfill $\Box$

The following lemma gives the estimate on the spatial gradients of the velocity, which is crucial for deriving the higher order estimates of the solution.
\begin{lemma}\label{lem35}
Under the condition \eqref{3.1}, it holds that for any $T\in[0,T^*)$,
\begin{equation}\label{5.1}
\sup_{0\leq t\leq T}\|\nabla\mathbf{u}\|_{L^2}^{2}
+\int_{0}^{T}\|\sqrt{\rho}\dot{\mathbf{u}}\|_{L^2}^{2}
 dt \leq C.
\end{equation}
\end{lemma}
{\it Proof.}
Multiplying \eqref{1.10}$_2$ by $\dot{\mathbf{u}}$ and integrating the resulting equation over $\mathbb{R}^2$ give rise to
\begin{align}\label{5.2}
\int\rho|\dot{\mathbf{u}}|^2dx
= -\int\dot{\mathbf{u}}\cdot\nabla Pdx+\mu\int\dot{\mathbf{u}}\cdot\Delta\mathbf{u}dx
+(\lambda+\mu)\int\dot{\mathbf{u}}\cdot\nabla\divv\mathbf{u}dx.
\end{align}
By \eqref{1.10}$_3$ and integrating by parts, we derive from \eqref{3.1} that
\begin{align}\label{5.3}
-\int\dot{\mathbf{u}}\cdot\nabla Pdx
& = \int \left[(\divv\mathbf{u})_tP
-(\mathbf{u}\cdot\nabla\mathbf{u})\cdot\nabla P\right]dx \nonumber\\
&=\frac{d}{dt}\int\divv\mathbf{u}Pdx+
\int\left[P(\divv\mathbf{u})^2-2\mu\divv\mathbf{u}|\mathfrak{D}\mathbf{u}|^2
-\lambda(\divv\mathbf{u})^3+P\partial_ju_i\partial_iu_j\right]dx \nonumber\\
&\leq \frac{d}{dt}\int\divv\mathbf{u}Pdx+
C\|\nabla\mathbf{u}\|_{L^2}^2+C\|\nabla\mathbf{u}\|_{L^3}^3.
\end{align}
It follows from integration by parts that
\begin{align}\label{5.4}
\mu\int\dot{\mathbf{u}}\cdot\Delta\mathbf{u}dx
& =\mu\int(\mathbf{u}_{t} +\mathbf{u}\cdot\nabla\mathbf{u})\cdot\Delta\mathbf{u}dx \notag \\
& = -\frac{\mu}{2}\frac{d}{dt}\|\nabla\mathbf{u}\|_{L^2}^{2}
-\mu\int\partial_{i}u_{j}\partial_{i}(u_{k}\partial_{k}u_{j})dx \notag \\
 &  \leq-\frac{\mu}{2}\frac{d}{dt}\|\nabla\mathbf{u}\|_{L^2}^{2}
 +C\|\nabla\mathbf{u} \|_{L^3}^{3}.
\end{align}
Similarly, one gets
\begin{align}\label{5.5}
(\lambda+\mu)\int\dot{\mathbf{u}}\cdot\nabla\divv\mathbf{u}dx
& = -\frac{\lambda+\mu}{2}\frac{d}{dt}\|\divv\mathbf{u}\|_{L^2}^{2}
-(\lambda+\mu)\int\divv\mathbf{u}\divv(\mathbf{u}\cdot\nabla\mathbf{u})dx \notag \\
 &  \leq-\frac{\lambda+\mu}{2}\frac{d}{dt}\|\divv\mathbf{u}\|_{L^2}^{2}
 +C\|\nabla\mathbf{u} \|_{L^3}^{3}.
\end{align}
Putting \eqref{5.3}--\eqref{5.5} into \eqref{5.2}, we obtain from \eqref{2.5} and \eqref{3.1} that
\begin{align}\label{5.6}
\Psi'(t)+\int\rho|\dot{\mathbf{u}}|^2dx
& \leq C\|\nabla\mathbf{u}\|_{L^2}^2+C\|\nabla\mathbf{u}\|_{L^3}^3 \notag \\
& \leq C\|\nabla\mathbf{u}\|_{L^2}^2+C\|\nabla\mathbf{u}\|_{L^3}^3
 \notag \\
& \leq C\|\nabla\mathbf{u}\|_{L^2}^2
+C\|\rho\dot{\mathbf{u}}\|_{L^2}\left(\|\nabla\mathbf{u}\|_{L^2}+1\right)^2
 \notag \\
& \leq \frac12\|\sqrt{\rho}\dot{\mathbf{u}}\|_{L^2}^2
+C(1+\|\nabla\mathbf{u}\|_{L^2}^2)\|\nabla\mathbf{u}\|_{L^2}^2+C,
\end{align}
where
\begin{equation*}
\Psi(t)\triangleq\frac{\mu}{2}\|\nabla\mathbf{u}\|_{L^2}^{2}
+\frac{\lambda+\mu}{2}\|\divv\mathbf{u}\|_{L^2}^{2}
-\int\divv\mathbf{u}Pdx
\end{equation*}
satisfies
\begin{equation}\label{5.7}
\frac{\mu}{2}\|\nabla\mathbf{u}\|_{L^2}^2-C
\leq \Psi(t)\leq\mu\|\nabla\mathbf{u}\|_{L^2}^2+C
\end{equation}
due to \eqref{3.3}.
Thus the desired \eqref{5.1} follows from \eqref{5.6}, \eqref{5.7}, \eqref{3.3},  and Gronwall's inequality. This completes the proof of Lemma \ref{lem35}.
\hfill $\Box$

Next, motivated by \cite{H1995}, we have the following estimates on the material derivatives of the velocity which are important for the higher order estimates of strong solutions.
\begin{lemma}\label{lem36}
Under the condition \eqref{3.1}, it holds that for any $T\in[0,T^*)$,
\begin{equation}\label{6.1}
\sup_{0\leq t\leq T}\|\sqrt{\rho}\dot{\mathbf{u}}\|_{L^2}^2
+\int_0^T\|\nabla\dot{\mathbf{u}}\|_{L^2}^2dt \leq C.
\end{equation}
\end{lemma}
{\it Proof.}
By the definition of $\dot{\mathbf{u}}$, we can rewrite $\eqref{1.10}_2$ as follows:
\begin{equation}\label{6.2}
\rho\dot{\mathbf{u}}+\nabla P
=\mu\Delta\mathbf{u}+(\lambda+\mu)\nabla\divv\mathbf{u}.
\end{equation}
Differentiating \eqref{6.2} with respect to $t$ and using \eqref{1.10}$_1$, we have
\begin{align}\label{6.3}
\rho\dot{\mathbf{u}}_t+\rho\mathbf{u}\cdot\nabla\dot{\mathbf{u}}
+\nabla P_t&=\mu\Delta\dot{\mathbf{u}}+(\lambda+\mu)\divv\dot{\mathbf{u}}
-\mu\Delta(\mathbf{u}\cdot\nabla\mathbf{u}) \nonumber \\
& \quad -(\lambda+\mu)\divv(\mathbf{u}\cdot\nabla\mathbf{u})
+\divv(\rho\dot{\mathbf{u}}\otimes\mathbf{u}).
\end{align}
Multiplying \eqref{6.3} by $\dot{\mathbf{u}}$ and integrating by parts over $\mathbb{R}^3$, we get
\begin{align}\label{6.4}
&\frac12\frac{d}{dt}\int \rho|\dot{\mathbf{u}}|^2\mbox{d}x+\mu\int|\nabla\dot{\mathbf{u}}|^2dx
+(\lambda+\mu)\int|\divv\dot{\mathbf{u}}|^2dx  \nonumber\\
&=\int \left(P_t\divv\dot{\mathbf{u}}+(\nabla P\otimes\mathbf{u}):\nabla\dot{\mathbf{u}}\right)dx
 +\mu\int [\divv(\Delta\mathbf{u}\otimes\mathbf{u})-\Delta(\mathbf{u}\cdot\nabla\mathbf{u})]\cdot \dot{\mathbf{u}}dx  \nonumber\\
& \quad+(\lambda+\mu)\int [(\nabla\divv\mathbf{u})\otimes\mathbf{u}
-\nabla\divv(\mathbf{u}\cdot\nabla\mathbf{u})]\cdot \dot{\mathbf{u}}dx\triangleq\sum_{i=1}^{3}J_i,
\end{align}
where $J_i$ can be bounded as follows.

It follows from $\eqref{1.10}_3$ that
\begin{align}\label{6.5}
J_1 & =\int \left(-\divv(P\mathbf{u})\divv\dot{\mathbf{u}}
-P\divv\mathbf{u}\divv\dot{\mathbf{u}}
+\mathcal{T}(\mathbf{u}):\nabla\mathbf{u}\divv\dot{\mathbf{u}}
+(\nabla P\otimes\mathbf{u}):\nabla\dot{\mathbf{u}}\right)dx   \nonumber\\
&=\int \left(P\mathbf{u}\nabla\divv\dot{\mathbf{u}}-P\divv\mathbf{u}\divv\dot{\mathbf{u}}
+\mathcal{T}(\mathbf{u}):\nabla\mathbf{u}\divv\dot{\mathbf{u}}
-P\nabla\mathbf{u}^\top:\nabla\dot{\mathbf{u}}
-P\mathbf{u}\nabla\divv\dot{\mathbf{u}}\right)dx \nonumber\\
&=\int \left(-P\divv\mathbf{u}\divv\dot{\mathbf{u}}
+\mathcal{T}(\mathbf{u}):\nabla\mathbf{u}\divv\dot{\mathbf{u}}
-P\nabla\mathbf{u}^\top:\nabla\dot{\mathbf{u}}\right)dx\nonumber\\
&\leq C\int \left(|\nabla\mathbf{u}||\nabla\dot{\mathbf{u}}|
+|\nabla\mathbf{u}|^2|\nabla\dot{\mathbf{u}}|\right)dx  \nonumber\\
&\leq C\left(\|\nabla\mathbf{u}\|_{L^2}+\|\nabla\mathbf{u}\|_{L^4}^2\right)
\|\nabla\dot{\mathbf{u}}\|_{L^2},
\end{align}
where $\mathcal{T}(\mathbf{u})=2\mu\mathfrak{D}(\mathbf{u})+\lambda\divv\mathbf{u}\mathbb{I}_3$.

For $J_2$ and $J_3$, notice that for all $1\leq i,j,k\leq 3,$ one has
\begin{align}
\partial_j(\partial_{kk}u_iu_j)-\partial_{kk}(u_j\partial_j u_i)&=\partial_k(\partial_ju_j\partial_ku_i)-\partial_k(\partial_ku_j\partial_ju_i)
-\partial_j(\partial_ku_j\partial_ku_i),\nonumber\\
\partial_j(\partial_{ik}u_ku_j)-\partial_{ij}(u_k\partial_ku_j)
&=\partial_i(\partial_ju_j\partial_ku_k)-\partial_i(\partial_ju_k\partial_ku_j)
-\partial_k(\partial_iu_k\partial_ju_j).\nonumber
\end{align}
So integrating by parts gives
\begin{align}
J_2&=\mu\int [\partial_k(\partial_ju_j\partial_ku_i)-\partial_k(\partial_ku_j\partial_ju_i)
-\partial_j(\partial_ku_j\partial_ku_i)]\dot{u_i}\mbox{d}x\nonumber\\
&\leq C\|\nabla\mathbf{u}\|_{L^4}^2\|\nabla\dot{\mathbf{u}}\|_{L^2},\label{6.7}\\
J_3&=(\lambda+\mu)\int [\partial_i(\partial_ju_j\partial_ku_k)-\partial_i(\partial_ju_k\partial_ku_j)
-\partial_k(\partial_iu_k\partial_ju_j)]\dot{u_i}\mbox{d}x\nonumber\\
&\leq C\|\nabla\mathbf{u}\|_{L^4}^2\|\nabla\dot{\mathbf{u}}\|_{L^2}.\label{6.8}
\end{align}
Inserting \eqref{6.5}--\eqref{6.8} into \eqref{6.4} and applying \eqref{5.1} lead to
\begin{align}\label{6.9}
& \frac12\frac{d}{dt}
\|\sqrt{\rho}\dot{\mathbf{u}}\|_{L^2}^2+\mu\|\nabla\dot{\mathbf{u}}\|_{L^2}^2
+(\lambda+\mu)\|\divv\dot{\mathbf{u}}\|_{L^2}^2  \nonumber\\
& \leq C(\|\nabla\mathbf{u}\|_{L^2}+\|\nabla\mathbf{u}\|_{L^4}^2)
\|\nabla\dot{\mathbf{u}}\|_{L^2}\nonumber\\
&\leq \frac{\mu}{2}\|\nabla\dot{\mathbf{u}}\|_{L^2}^2 +C(\mu)\left(\|\nabla\mathbf{u}\|_{L^4}^4+1\right),
\end{align}
which implies
\begin{align}\label{6.16}
\frac{d}{dt}\|\sqrt{\rho}\dot{\mathbf{u}}\|_{L^2}^2
+\mu\|\nabla\dot{\mathbf{u}}\|_{L^2}^2
\leq C\|\nabla\mathbf{u}\|_{L^4}^4+C.
\end{align}
By virtue of \eqref{2.5}, \eqref{3.2}, \eqref{3.3}, and \eqref{5.1}, one has
\begin{equation}\label{6.17}
\|\nabla\mathbf{u}\|_{L^4}^4
\leq C\|\rho\dot{\mathbf{u}}\|_{L^2}^2
\left(\|\nabla\mathbf{u}\|_{L^2}+\|P\|_{L^2}\right)^2
+C\|P\|_{L^4}^4
\leq C\|\sqrt{\rho}\dot{\mathbf{u}}\|_{L^2}^2
+C.
\end{equation}
Consequently, we obtain the desired \eqref{6.1} from \eqref{6.16}, \eqref{6.17}, and Gronwall's inequality. This completes the proof of Lemma \ref{lem36}.
\hfill $\Box$

Inspired by \cite{LX2013,LSZ2015}, we have
the following spatial weighted estimate on the density, which plays an important role in deriving the bounds on the higher order derivatives of the solutions $(\rho,\mathbf{u},P)$.
\begin{lemma}\label{lem36}
Under the condition \eqref{3.1}, it holds that for any $T\in[0,T^*)$,
\begin{equation}\label{06.1}
\sup_{0\leq t\leq T}\|\rho\bar{x}^{a}\|_{L^{1}}\leq C.
\end{equation}
\end{lemma}

{\it Proof.}
First, for $N>1,$ let $\varphi_N\in C^\infty_0(\mathbb{R}^2)$  satisfy
\begin{equation} \label{vp1}
0\le \varphi_N \le 1, \quad  \varphi_N(x)
=\begin{cases} 1,~~~~ |x|\le N/2,\\
0,~~~~ |x|\ge N,\end{cases}
\quad |\nabla \varphi_N|\le C N^{-1}.
\end{equation}
It follows from \eqref{1.10}$_1$ that
\begin{align}\label{oo0}
\frac{d}{dt}\int \rho \varphi_{N} dx &=\int \rho\mathbf{u}
\cdot\nabla \varphi_{N} dx \notag \\
&\ge - C N^{-1}\left(\int\rho dx\right)^{1/2}
\left(\int\rho |u|^2dx\right)^{1/2}\ge - \tilde{C} N^{-1},
\end{align}
where in the last inequality one has used \eqref{3.2} and \eqref{3.3}.
Integrating \eqref{oo0} and choosing $N=N_1\triangleq2N_0+4\tilde CT$, we obtain after using \eqref{1.9} that
\begin{align}\label{p1}
\inf\limits_{0\le t\le T}\int_{B_{N_1}} \rho dx&\ge \inf\limits_{0\le t\le T}\int \rho \varphi_{N_1} dx \notag \\
&\ge \int \rho_0 \varphi_{N_1} dx-\tilde{C}N_1^{-1}T \notag \\
&\ge \int_{B_{N_0}} \rho_0 dx-\frac{\tilde{C}T}{2N_0+4\tilde{C} T} \notag \\
&\ge 1/4.
\end{align}
Hence, it follows from \eqref{p1}, \eqref{3.2}, \eqref{22}, \eqref{3.3}, and \eqref{5.1} that for any $\eta\in(0,1]$ and any $s>2$,
\begin{equation}\label{06.2}
\|\mathbf{u}\bar{x}^{-\eta}\|_{L^{s/\eta}}\leq C\left(\|\sqrt{\rho}\mathbf{u}\|_{L^2}+\|\nabla\mathbf{u}\|_{L^2}\right)\le C.
\end{equation}
Multiplying \eqref{1.10}$_1$ by $\bar{x}^{a}$ and integrating the resulting equality by parts over $\mathbb{R}^2$ yield that
\begin{equation*}
\begin{split}
\frac{d}{dt}\int\rho\bar{x}^{a}dx & \leq C\int\rho|\mathbf{u} |\bar{x}^{a-1}\log^{2}(e+|x|^2)dx\\
 &  \leq C\|\rho\bar{x}^{a-1+\frac{8}{8+a}}\|_{L^{\frac{8+a}{7+a}}}\|\mathbf{u} \bar{x}^{-\frac{4}{8+a}}\|_{L^{8+a}} \\
 &  \leq C\int\rho\bar{x}^{a}dx+C,
\end{split}
\end{equation*}
which along with Gronwall's inequality gives \eqref{06.1} and finishes the proof of Lemma \ref{lem36}.    \hfill $\Box$

The following lemma will treat the higher order derivatives of the solutions which are needed to guarantee the extension of local strong solution to be a global one.
\begin{lemma}\label{lem37}
Under the condition \eqref{3.1}, and let $q>2$ be as in Theorem \ref{thm1.1}, then it holds that for any $T\in[0,T^*)$,
\begin{equation}\label{7.1}
\sup_{0\leq t\leq T}\left(\|(\rho,P)\|_{H^1\cap W^{1,q}}+\|\nabla\mathbf{u}\|_{H^1}\right)
+\int_{0}^T\|\nabla^2\mathbf{u}\|_{L^q}^2dt\leq C.
\end{equation}
\end{lemma}
{\it Proof.}
First, it follows from the mass equation \eqref{1.10}$_1$ that $\nabla\rho$ satisfies for any $r\in[2,q]$,
\begin{align}\label{7.2}
\frac{d}{dt}\|\nabla\rho\|_{L^r}
& \leq C(r)(1+\|\nabla\mathbf{u} \|_{L^\infty})\|\nabla\rho\|_{L^r}
+C(r)\|\nabla^2\mathbf{u}\|_{L^r} \nonumber \\
& \leq C(1+\|\nabla\mathbf{u} \|_{L^\infty})\|\nabla\rho\|_{L^r}
+C(\|\rho\dot{\mathbf{u}}\|_{L^r}+\|\nabla P\|_{L^r})
\end{align}
due to
\begin{equation}\label{7.3}
\|\nabla^2\mathbf{u}\|_{L^r}
\leq C(\|\rho\dot{\mathbf{u}}\|_{L^r}+\|\nabla P\|_{L^r}),
\end{equation}
which follows from the standard $L^r$-estimate for the following elliptic system
\begin{align*}
\begin{cases}
\mu\Delta\mathbf{u}+(\lambda+\mu)\nabla\divv\mathbf{u}=\rho\dot{\mathbf{u}}+\nabla P,\ \ x\in\mathbb{R}^2,\\
\mathbf{u}\rightarrow\mathbf{0},\ \text{as}\ |x|\rightarrow\infty.
\end{cases}
\end{align*}
Similarly, one deduces from \eqref{1.10}$_3$ that $\nabla P$ satisfies for any $r\in[2,q]$,
\begin{align}\label{7.4}
\frac{d}{dt}\|\nabla P\|_{L^r} & \leq C(r)(1+\|\nabla\mathbf{u}\|_{L^{\infty}})(\|\nabla P\|_{L^r}+\|\nabla^2\mathbf{u}\|_{L^r}) \notag \\
& \leq C(1+\|\nabla\mathbf{u}\|_{L^{\infty}})(\|\rho\dot{\mathbf{u}}\|_{L^r}+\|\nabla P\|_{L^r}).
\end{align}

Next, one gets from \eqref{3.14}, Gagliardo-Nirenberg inequality, \eqref{3.1}, \eqref{2.3}, \eqref{6.1}, and \eqref{3.2} that
\begin{align}\label{7.5}
\|\divv\mathbf{u}\|_{L^\infty}+\|\omega\|_{L^\infty}
& \leq C\|P\|_{L^\infty}+C\|F\|_{L^\infty}+\|\omega\|_{L^\infty} \nonumber \\
& \leq C(q)+C(q)\|\nabla F\|_{L^2}^{\frac{q-2}{2(q-1)}}\|\nabla F\|_{L^q}^{\frac{q}{2(q-1)}}
+C(q)\|\nabla\omega\|_{L^2}^{\frac{q-2}{2(q-1)}}\|\nabla\omega\|_{L^q}^{\frac{q}{2(q-1)}} \nonumber \\
& \leq C+C\|\rho\dot{\mathbf{u}}\|_{L^q}^{\frac{q}{2(q-1)}},
\end{align}
which together with Lemma \ref{lem25}, \eqref{7.3}, and \eqref{5.1} yields that
\begin{align}\label{7.6}
\|\nabla\mathbf{u}\|_{L^\infty}
& \leq C\left(\|\divv\mathbf{u}\|_{L^\infty}+\|\omega\|_{L^\infty}\right)
\log(e+\|\nabla^2\mathbf{u}\|_{L^q})+C\|\nabla\mathbf{u}\|_{L^2}+C \nonumber \\
& \leq C\left(1+\|\rho\dot{\mathbf{u}}\|_{L^q}^{\frac{q}{2(q-1)}}\right)
\log\left(e+\|\rho\dot{\mathbf{u}}\|_{L^q}+\|\nabla P\|_{L^q}\right)+C.
\end{align}
It follows from \eqref{p1}, \eqref{3.1}, \eqref{22}, and \eqref{06.1} that for any $\eta\in(0,1]$ and any $s>2$,
\begin{align}\label{7.7}
\|\rho^\eta v\|_{L^{\frac{s}{\eta}}}
 &  \leq C\|\rho^\eta\bar x^{\frac{3\eta a}{4s}}\|_{L^{\frac{4s}{3\eta}}}
\|v\bar  x^{-\frac{3\eta a}{4s}}\|_{L^{\frac{4s}{\eta}}} \notag \\
 &  \leq C\|\rho\|_{L^\infty}^{\frac{(4s-3)\eta}{4s}}\|\rho\bar x^a\|_{L^1}^{\frac{3\eta}{4s}}\left( \|\sqrt{\rho} v\|_{L^2}+\|\nabla v\|_{L^2}\right) \notag \\
 &  \leq C\left(\|\sqrt{\rho}v\|_{L^2}+\|\nabla v\|_{L^2}\right),
\end{align}
which along with H{\"o}lder's inequality, \eqref{6.1}, and \eqref{3.2} shows that
\begin{align}\label{7.8}
\|\rho\dot{\mathbf{u}}\|_{L^q}
 &  \leq C\|\rho\dot{\mathbf{u}}\|_{L^2}^{\frac{2(q-1)}{q^2-2}}
\|\rho\dot{\mathbf{u}}\|_{L^{q^2}}^{\frac{q(q-2)}{q^2-2}} \notag \\
 &  \leq C\|\rho\dot{\mathbf{u}}\|_{L^2}^{\frac{2(q-1)}{q^2-2}}
\left(\|\sqrt{\rho}\dot{\mathbf{u}}\|_{L^2}+\|\nabla\dot{\mathbf{u} }\|_{L^2}\right)^{\frac{q(q-2)}{q^2-2}} \notag \\
 &  \leq C\left(1+\|\nabla\dot{\mathbf{u}}\|_{L^2}^{\frac{q(q-2)}{q^2-2}}\right).
\end{align}
Then we derive from \eqref{7.6} and \eqref{7.8} that
\begin{align}\label{7.9}
\|\nabla\mathbf{u}\|_{L^\infty}
\leq C\left(1+\|\nabla\dot{\mathbf{u}}\|_{L^2}\right)
\log\left(e+\|\nabla\dot{\mathbf{u}}\|_{L^2}+\|\nabla P\|_{L^q}\right)+C
\end{align}
due to $\frac{q(q^2-2q)}{(2q-2)(q^2-2)},\ \frac{q^2-2q}{q^2-2}\in(0,1)$.
Consequently, substituting \eqref{7.8} and \eqref{7.9}  into \eqref{7.2} and \eqref{7.4}, we get after choosing $r=q$ that
\begin{equation}\label{7.11}
f'(t)\leq  Cg(t)f(t)\log{f(t)}+Cg(t)f(t)+Cg(t),
\end{equation}
where
\begin{align*}
f(t)&\triangleq e+\|\nabla\rho\|_{L^{ q}}+\|\nabla P\|_{L^{ q}},\\
g(t)&\triangleq (1+\|\nabla\dot{\mathbf{u}}\|_{L^2})
\log(e+\|\nabla\dot{\mathbf{u}}\|_{L^2}).
\end{align*}
This yields
\begin{equation}\label{7.12}
(\log f(t))'\leq Cg(t)+Cg(t)\log f(t)
\end{equation}
due to $f(t)>1.$
Thus it follows from  \eqref{7.12}, \eqref{6.1}, and  Gronwall's inequality that
\begin{equation}\label{7.13}
\sup_{0\leq t\leq T}\|(\nabla\rho,\nabla P)\|_{L^{q}}\leq C,
\end{equation}
which, combined with \eqref{7.6}, \eqref{7.8}, and \eqref{6.1} gives that
\begin{equation}\label{7.14}
\int_{0}^{T}\|\nabla\mathbf{u}\|^2_{L^\infty}dt\leq C.
\end{equation}
Taking $r=2$ in \eqref{7.2} and \eqref{7.4}, one gets from \eqref{7.14}, \eqref{7.7}, \eqref{6.1}, and Gronwall's inequality  that
\begin{equation}\label{7.15}
\sup_{0\leq t\leq T}\|(\nabla\rho,\nabla P)\|_{L^2}\leq C,
\end{equation}
which together with \eqref{7.3}, \eqref{6.1}, and \eqref{3.2} yields that
\begin{equation}\label{7.16}
\sup_{0\leq t\leq T}\|\nabla^2\mathbf{u}\|_{L^2}\leq C.
\end{equation}
Taking $r=q$ in \eqref{7.3} and using \eqref{5.1} and \eqref{7.13} show that
\begin{equation}\label{7.17}
\int_{0}^T\|\nabla^2\mathbf{u}\|_{L^q}^2dt\leq C.
\end{equation}
Thus the desired \eqref{7.1} follows from \eqref{7.17}, \eqref{7.13}, \eqref{7.15}, \eqref{7.16}, \eqref{5.1}, \eqref{3.2}, and \eqref{3.3}. The proof of Lemma \ref{lem37} is finished.  \hfill $\Box$

%which along with \eqref{4.5} in particular implies
%\begin{equation} \label{nu1}
%\int_0^T\|\nabla\mathbf{u}\|_{L^\infty}dt\leq C.
%\end{equation}
%Thus, applying Gronwall's inequality to \eqref{t} gives
%\begin{equation}\label{4.8}
%\sup_{t\in[0,T]}\|\nabla\rho\|_{L^2\cap L^q}\leq C.
%\end{equation}

\begin{lemma}\label{lem38}
Under the condition \eqref{3.1}, it holds that for any $T\in[0,T^*)$,
\begin{equation}\label{8.1}
\sup_{0\leq t\leq T}\|\rho\bar{x}^{a}\|_{L^1\cap H^{1}\cap W^{1,q}}\leq C.
\end{equation}
\end{lemma}

{\it Proof.}
One derives from \eqref{1.10}$_1$ that $\rho\bar{x}^a$ satisfies
\begin{equation}\label{8.2}
\partial_{t}(\rho\bar{x}^a)+ \mathbf{u}\cdot\nabla(\rho\bar{x}^a)
-a\rho\bar{x}^{a}\mathbf{u}\cdot\nabla\log\bar{x}+\rho\bar{x}^a\divv\mathbf{u}=0.
\end{equation}
Taking the $x_i$-derivative on the both side of \eqref{8.2} gives
\begin{align}\label{8.3}
0=  &  \partial_{t}\partial_{i}(\rho\bar{x}^a)+\mathbf{u} \cdot\nabla\partial_{i}(\rho\bar{x}^a)
+\partial_{i}\mathbf{u}\cdot\nabla(\rho\bar{x}^a)
-a\partial_{i}(\rho\bar{x}^{a})\mathbf{u}\cdot\nabla\log\bar{x} \notag \\
 &  -a\rho\bar{x}^{a}\partial_{i}\mathbf{u}\cdot\nabla\log\bar{x}
-a\rho\bar{x}^{a} \mathbf{u}\cdot\partial_{i}\nabla\log\bar{x}
+\partial_{i}(\rho\bar{x}^a\divv\mathbf{u}).
\end{align}
For any $r\in[2,q]$, multiplying \eqref{8.3} by $|\nabla (\rho\bar{x}^a)|^{r-2}\partial_{i}(\rho\bar{x}^a)$ and integrating the resulting equality over $\mathbb{R}^2$, we obtain from integrating by parts, \eqref{06.2}, \eqref{06.1}, \eqref{2.5}, \eqref{6.1}, and \eqref{7.1} that
\begin{align}\label{8.4}
\frac{d}{dt}\|\nabla(\rho\bar{x}^a)\|_{L^r} \leq  &
C\left(1+\|\nabla\mathbf{u}\|_{L^\infty}+\|\mathbf{u}
\cdot\nabla\log\bar{x}\|_{L^\infty}\right)
\|\nabla(\rho\bar{x}^a)\|_{L^r} \notag \\
&  + C\|\rho\bar{x}^a\|_{L^\infty}\left(\||\nabla\mathbf{u} ||\nabla\log\bar{x}|\|_{L^r}
+\||\mathbf{u}||\nabla^{2}\log\bar{x}|\|_{L^r}+\|\nabla^2\mathbf{u}\|_{L^r}\right) \notag \\
\leq & C\left(1+\|\nabla\mathbf{u}\|_{W^{1,q}}\right)\|\nabla(\rho\bar{x}^a)\|_{L^r}  \notag \\
& +C\|\rho\bar{x}^a\|_{L^\infty}\left(\|\nabla\mathbf{u}\|_{L^r}
+\|\mathbf{u}\bar{x}^{-\frac{2}{5}}\|_{L^{4r}}\|\bar{x}^{-\frac{3}{2}}
\|_{L^{\frac{4r}{3}}}+\|\nabla^2\mathbf{u}\|_{L^r}\right) \notag \\
\leq & C\left(1+\|\nabla^2\mathbf{u}\|_{L^r}+\|\nabla\mathbf{u}\|_{W^{1,q}}\right)
\left(1+\|\nabla(\rho\bar{x}^a)\|_{L^r}+\|\nabla(\rho\bar{x}^a)\|_{L^q}\right)
\notag \\
\leq & C\left(1+\|\nabla^2\mathbf{u}\|_{L^r}+\|\nabla^2\mathbf{u}\|_{L^q}\right)
\left(1+\|\nabla(\rho\bar{x}^a)\|_{L^r}+\|\nabla(\rho\bar{x}^a)\|_{L^q}\right).
\end{align}
Choosing $r=q$ in \eqref{8.4}, together with \eqref{7.1} and Gronwall's inequality indicates that
\begin{equation}\label{8.5}
\sup_{t\in[0,T]}\|\nabla(\rho\bar{x}^a)\|_{L^q}\leq C.
\end{equation}
Setting $r=2$ in \eqref{8.4}, we deduce from \eqref{7.1} and \eqref{8.5} that
\begin{equation*}
\sup_{t\in[0,T]}\|\nabla(\rho\bar{x}^a)\|_{L^2}\leq C.
\end{equation*}
This combined with \eqref{8.5} and \eqref{06.1} gives \eqref{8.1} and completes the proof of Lemma \ref{lem38}.   \hfill $\Box$

With Lemmas \ref{lem31}--\ref{lem38} at hand, we are now in a position to prove Theorem \ref{thm1.1}.

\textbf{Proof of Theorem \ref{thm1.1}.}
We argue by contradiction. Suppose that \eqref{B} were false, that is, \eqref{3.1} holds. Note that the general constant $C$ in Lemmas \ref{lem31}--\ref{lem38} is independent of $t<T^{*}$, that is, all the a priori estimates obtained in Lemmas \ref{lem31}--\ref{lem38} are uniformly bounded for any $t<T^{*}$. Hence, the function
\begin{equation*}
(\rho,\mathbf{u},P)(x,T^{*})
\triangleq\lim_{t\rightarrow T^{*}}(\rho,\mathbf{u},P)(x,t)
\end{equation*}
satisfy the initial condition \eqref{2.02} at $t=T^{*}$.

Furthermore, standard arguments yield that $\rho\dot{\mathbf{u}}\in C([0,T];L^2)$, which
implies $$ \rho\dot{\mathbf{u}}(x,T^\ast)=\lim_{t\rightarrow
T^\ast}\rho\dot{\mathbf{u}}\in L^2. $$
Hence, $$-\mu\Delta{\mathbf{u}}-(\lambda+\mu)\nabla\mbox{div}\mathbf{u}+\nabla P
|_{t=T^\ast}=\sqrt{\rho}(x,T^\ast)g(x)
$$ with $$g(x)\triangleq
\begin{cases}
\rho^{-1/2}(x,T^\ast)(\rho\dot{\mathbf{u}})(x,T^\ast),&
\mbox{for}~~x\in\{x|\rho(x,T^\ast)>0\},\\
0,&\mbox{for}~~x\in\{x|\rho(x,T^\ast)=0\},
\end{cases}
$$
satisfying $g\in L^2$ due to \eqref{7.1}.
Therefore, one can take $(\rho,\mathbf{u},P)(x,T^\ast)$ as
the initial data and extend the local
strong solution beyond $T^\ast$. This contradicts the assumption on
$T^{\ast}$.

Thus we finish the proof of Theorem \ref{thm1.1}.
\hfill $\Box$

\end{document}